\documentclass[twoside]{aiml22}

\usepackage{aiml22macro}

\usepackage{graphicx}
\usepackage{amsmath}
\usepackage{amssymb}

\def\lastname{Almeida and Bezhanishvili}

\begin{document}

\begin{frontmatter}
  \title{Coalgebraic Semantics for  Intuitionistic \\ Modal Logic}
  \author{Rodrigo Nicolau Almeida}\footnote{r.dacruzsilvapinadealmeida@uva.nl.}
  \address{Institute for Logic Language and Computation (ILLC) - University of Amsterdam \\ Science Park 107 \\ 1098 XG Amsterdam}
  \author{Nick Bezhanishvili}
  \address{Institute for Logic Language and Computation (ILLC) - University of Amsterdam \\ Science Park 107 \\ 1098 XG Amsterdam}

  \begin{abstract}
We give a new coalgebraic semantics for intuitionistic modal logic with $\Box$. In particular, we provide a colagebraic representation of  intuitionistic descriptive modal frames and  of intuitonistic modal Kripke frames based on 
image-finite posets.
This gives a solution to a problem in the area of coalgebaic logic  for these classes of frames,  raised  explicitly by Litak (2014) and de Groot and Pattinson (2020). Our key technical tool is a recent generalization of a construction by Ghilardi, in the form of a right adjoint to the inclusion of the category of Esakia spaces in the category of Priestley spaces. As an application of these results, we study bisimulations of intuitionistic modal frames, describe dual spaces of free modal Heyting algebras, and provide a path towards a theory of coalgebraic intuitionistic logics.  \end{abstract}

 \end{frontmatter}

\section{Introduction}

Coalgebraic semantics of classical modal and positive modal logics have been thoroughly investigated \cite{Venema2014modallogicandvietoris,Palmigiano2004,Venema2007,Kupke2003}. 
In the case of classical modal logic this is done via the powerset coalgebras in the context of Kripke frames and by the Vietoris coalgebras for descriptive frames (see \cite{Kupke2003,Venema2014modallogicandvietoris}). For positive modal logic this is obtained via the convex set functor and its topological analogue \cite{Balan2015,BezhanishviliK07,Palmigiano2004,Bezhanishvili2023remarksonhyperspaces}. 

 Intuitonistic modal logics, on the other hand, have so far escaped this kind of analysis. Such logics are quite varied, stemming from several distinct proposals over what is the appropriate notion of a ``constructive" version of modal logic, see, e.g.,  \cite{simpson1994proof,FischerServi1977,Servi1980,Esakia2006,Wolter1999}. They all have well-known algebraic semantics in the form of modal Heyting algebras,  a Kripke-style semantics \cite{Esakia2006,Litak2014} via so-called ``Intuitionistic modal frames" \cite{Wolter1999} and an order-topological semantics in terms of intuitionistic descriptive modal frames (see e.g., \cite{Wolter1999}, where these are called $\Box$-frames, $\dia$-frames, etc, depending on the specific signature used). However, the question of how to represent such frames as coalgebras for an appropriate functor on the category of posets with p-morphisms or respecitively on the category of Esakia spaces and continuous p-morphisms (which is dual to the category of Heyting algebras) was open in the research community for a long time. It was raised explicitly by Litak \cite{Litak2014} and by de Groot and Pattinson\footnote{We point out that de Groot and Pattinson \cite{deGroot2020} also use the so-called dialgebraic representations in their work.} \cite{deGroot2020}, where the existence of such a representation was left as an open question. 

In this paper we resolve the problem of defining coalgebraic semantics for intuitionistic modal logic, albeit for particular classes of intuitionistic modal frames. More concretely, we will represent intuitionistic descriptive modal frames and intuitionistic modal Kripke frames based on 
image-finite posets as particular coalgebras. Intuitionistic descriptive modal frames can be seen as Esakia spaces equipped with a modal relation. 
We will represent these as particular coalgebras on the category of Esakia spaces.
For this we will 
be 
%
%This paper seeks to address the question of how to develop a coalgebraic semantics for several intuitionistic modal logics. 
%In this paper we solve this problem, 
relying on a recent generalization \cite{almeidanoteoncolimitsofHeytingalgebras} of a construction by Ghilardi \cite{ghilardifreeheyting}. As our main results show, descriptive intuitionistic modal frames can be seen as coalgebras for an appropriate composition of a variant of the classical upwards-Vietoris set functor and a right adjoint to the inclusion of the category of Esakia spaces with continuous p-morphisms in the category of Priestley spaces. Similar results are likewise obtained for the categories of image-finite posets with p-morphisms and posets with monotone maps, with appropriate modifications, giving us a coalgebraic representation of all image-finite intuitionistic modal frames. We observe,  however, that this does not provide a representation of all intuitionistic modal frames, since the usage of image-finiteness is quite crucial. 
%a situation which we briefly discuss.

We provide three applications of these constructions: (1) a colagebraic notion of a bisimulation for intuitionistic descrptive and image-finite Kripke frames, (2) a concrete description of the dual space of the free modal Heyting algebra on finitely many generators, and (3) a new definition for coalgebraic intuitionistic logic laying a path towards a theory of such logics following the line of research suggested at the end of \cite{deGroot2020}.

The outline of the paper is as follows: in Section \ref{Preliminaries} we recall the necessary preliminaries from the theory of coalgebras over $\mathbf{Set}$ and $\mathbf{Stone}$, as well as the relevant versions of intuitionistic modal logic and its semantics. In Section \ref{Generalizing Ghilardi} we review the key aspects of the construction from \cite{almeidanoteoncolimitsofHeytingalgebras} we will need for our purpose. We provide our main results showcasing the coalgebraic semantics, and its equivalence with the classical descriptive frame semantics, in Section \ref{Main Results}. We present our applications in Section \ref{Applications}. We conclude in Section \ref{Conclusions} by pointing to further research directions in the study of these representations.

\section{Preliminaries}\label{Preliminaries}

\subsection{Duality Theory and Coalgebra}

We assume throughout that the reader is familiar with Stone and Priestley duality (as presented, e.g. in \cite{Davey2002-lr}). 

\begin{definition}
    Let $(X,\leq,\tau)$ be an ordered-topological space. We say that $X$ is a \textit{Priestley space} if $(X,\tau)$ is compact, and it satisfies the \textit{Priestley Separatiom Axiom}: if $x\nleq y$, there is a clopen upset $U$ such that $x\in U$ and $y\notin U$.
\end{definition}

\begin{definition}
    Let $(X,\leq),(Y,\leq)$ be two posets, and $f:X\to Y$ be a map between them. The map $f$ is said to be \textit{monotone} if whenever $x,y\in X$, and $x\leq y$ then $f(x)\leq f(y)$. We say that $f$ is a \textit{p-morphism} if it is monotone and in addition, whenever $x\in X$, $y\in Y$ and $f(x)\leq y$, then there is some $x'$ such that $x\leq x'$ and $f(x')=y$.
\end{definition}

We say that a map $f:X\to Y$ between Priestley spaces is a \textit{Priestley morphism} if it is continuous and monotone. We denote by $\mathbf{Pries}$ the category of Priestley spaces, with Priestley morphisms. We will in particular need the restriction of such a duality to \textit{Esakia duality}, which we now recall:

\begin{definition}
    Let $(X,\leq)$ be a Priestley space. We say that $X$ is an \textit{Esakia space} if whenever $U$ is a clopen set in $X$, then ${\downarrow}U$ is clopen as well.

    Given Esakia spaces $X,Y$, a morphism $f:X\to Y$ is said to be an \textit{Esakia morphism} if it is a Priestley morphism and a p-morphism between the underlying posets. We denote by $\mathbf{Esa}$ the category of Esakia spaces with Esakia morphisms.
\end{definition}

It is well-known (see e.g. \cite{Esakiach2019HeyAlg}) that the category $\mathbf{Esa}$ is dual to the category of Heyting algebras and Heyting algebra homomorphisms. We also recall J\'onnson-Tarski duality (see e.g. \cite{Blackburn2002-fd,Chagrov1997-cr}):

\begin{definition}
    Let $(X,R)$ be a Stone space where $R\subseteq X\times X$. We say that $X$ is a \textit{modal space} if:
    \begin{enumerate}
        \item For each $x\in X$, the set $R[x]$ is closed;
        \item For each clopen set $U$, $R^{-1}[U]$ is clopen.
    \end{enumerate}
    Given a modal space $(X,R)$, let $\mathfrak{X}=(X,R,\mathcal{A})$ be the triple where $\mathcal{A}=\mathsf{Clop}(X)$ is the set of clopen subsets. We call $\mathfrak{X}$ a \textit{descriptive general frame over $X$}. Given $f:X\to Y$ a map between descriptive general frames, we say that it is a descriptive morphism if $f$ is continuous, and whenever $x\in X$ and $y\in Y$, and $f(x)Ry$, then there is some $x'\in X$ such that $xRx'$ and $f(x')=y$. We denote by $\mathbf{DG}$ the category of descriptive general frames with descriptive morphisms.\footnote{We note that this definition of descriptive general frames is equivalent to the usual one that can be found in e.g. \cite{Blackburn2002-fd,Chagrov1997-cr}.}
\end{definition}

In parallel with these categories, we will work throughout with some categories of posets. Recall that a poset $(P,\leq)$ is said to be \textit{image-finite} if for each $x\in P$, ${\uparrow}x=\{y : x\leq y\}$ is finite. We will work with the following:
\begin{enumerate}
    \item $\mathbf{Pos}$, the category of posets with monotone maps;
    \item The (non-full) subcategory $\mathbf{Pos}_{p}$ of $\mathbf{Pos}$ where we restrict maps to p-morphisms;
    \item The subcategory $\mathbf{ImFinPos}_{p}$ of \textit{image-finite} posets with p-morphisms.
\end{enumerate}
As we will note below (see Section \ref{The functor PG} for further discussion), the categories $\mathbf{Pos}$ and $\mathbf{ImFinPos}_{p}$ play a similar role to the categories $\mathbf{Pries}$ and $\mathbf{Esa}$, respectively, in a discrete setting.

\begin{definition}
    Given a category $\mathbb{C}$, and an endofunctor $F:\mathbb{C}\to \mathbb{C}$, a pair $(A,f)$ of an object $A$ and a morphism $f:A\to F(A)$ is called an $F$-\textit{coalgebra} (or just a coalgebra, if the relevant $F$ is clear from context).

    Given two $F$-coalgebras $(A,f_{A})$ and $(B,f_{B})$, we say that a morphism $h:A\to B$ is a \textit{coalgebra morphism} between $(A,f_{A})$ and $(B,f_{B})$ if it makes the following diagram commute:

\begin{figure}[h]
        \centering
\begin{tikzcd}
A \arrow[d, "f_{A}"'] \arrow[r, "h"] & B \arrow[d, "f_{B}"] \\
F(A) \arrow[r, "F(h)"']              & F(B)                
\end{tikzcd}        \caption{Coalgebra morphism compatibility}
        \label{fig:coalgebramorphismcompatibility}
    \end{figure}
Such morphisms compose in the obvious way. We write $\mathbf{CoAlg}(F)$ for the category of $F$-coalgebras and coalgebra morphisms.    
\end{definition}

We refer the reader to \cite{Venema2007} for all of the facts about coalgebra we will assume here. If one considers in particular the category $\mathbf{Set}$ of sets and functions, the powerset functor $\mathcal{P}$ admits a particularly transparent description of its coalgebras: they are precisely the Kripke frames. Tracing its origins in the work of Esakia \cite{esakiatopologicalkripkemodels} (see also \cite{Venema2014modallogicandvietoris} for an in-depth discussion), it has been realised that also \textit{descriptive general frames} can be represented as coalgebras of a specific endofunctor on the category $\mathbf{Stone}$ of Stone spaces and continuous functions:

\begin{definition}
    Let $X$ be a Stone space. Let $V(X)$ be the set of closed subset of $X$. We give this set a topology consisting of the ``hit-and-miss" topology, i.e., by giving it the topology determined by the subbasis consisting of
    \begin{equation*}
        [U]=\{C\in V(X) : C\subseteq U\} \text{ and } \langle V\rangle =\{C\in V(X) : C\cap V\neq \emptyset\}
    \end{equation*}
    where $U,V$ range over clopen subsets of $X$. We call this space the \textit{Vietoris hyperspace} of $X$.    Moreover, given a continuous function $f:X\to Y$, we define $V(f)$ to be the direct image of $f$.
\end{definition}

Then we have the following (see e.g. \cite{Kupke2003,Venema2014modallogicandvietoris}):

\begin{proposition}
    The assignment $V$ as above defines an endofunctor on $\mathbf{Stone}$. Moreover, for each $X$ a Stone space, the categories $\mathbf{CoAlg}(V)$ and $\mathbf{DG}$ are equivalent.
\end{proposition}

\subsection{Intuitionistic Modal Logic}

We consider the language of $\mathsf{IPC}_{\Box}$ of intuitionistic logic with an additional unary operation $\Box$. Throughout this section we refer to \textit{intuitionistic modal logic}, $\mathsf{IPC}_{\Box}$, as the logic axiomatised by:
\begin{enumerate}
    \item The axioms from $\mathsf{IPC}$;
    \item $\Box(\phi\wedge \psi)\leftrightarrow \Box \phi\wedge \Box \psi$;
    \item $\Box\top\leftrightarrow \top$.
\end{enumerate}

Using the usual completeness methods, one can show that this logic is complete with respect to some Kripke-style semantics. The relevant semantics for such a logic is given over modal intuitionistic frames (sometimes called $\Box$-frames, to distinguish them from the semantics of richer intuitionistic modal logics):

\begin{definition}
    Let $(X,\leq,R)$ be a triple where $(X,\leq)$ is a partial order, $R\subseteq X\times X$. We say this is a \textit{modal intuitionistic frame} if it satisfies
    \begin{equation*}
        {R}{=}{\leq} \circ {R} \circ {\leq}.
    \end{equation*}
    Given two modal intuitionistic frames $(X,\leq,R)$ and $(Y,\leq,R)$ we say that a map $f:X\to Y$ is a \textit{modal p-morphism} if $f$ is a p-morphism with respect to $\leq$ and $R$.
    Let $\mathbf{ImFinK}$ be the category of image-finite modal intuitionistic frames with modal p-morphisms.
\end{definition}

Models are constructed by taking valuations in $\mathsf{Up}(X)$, the set of upsets of the poset. The semantics of the intuitionistic connectives is kept the same, whilst the clause for the $\Box$-operator is similar to classical modal logic: a model $\bbM,x\sat \Box\phi$ if and only if whenever $y\in R[x]$ then $\bbM,y\sat \phi$.

Note that if $X$ is a poset, then $\mathsf{Up}(X)$  carries a richer structure: ordering upwards closed subsets using \textit{reverse inclusion}, we obtain that $(\mathsf{Up}(X),\supseteq)$ is a poset. And such an assignment is in fact an endofunctor on the category $\mathbf{Pos}$, sending monotone maps to their direct image.

Now, if $f:X\to \mathsf{Up}(X)$ is a coalgebra for this functor, then we can think of this map as picking for each $x\in X$ a set of modal successors; and dually, given a modal intuitionistic frame $(X,\leq,R)$, the map $R_{X}:X\to \mathsf{Up}(X)$ will provide a coalgebra. Hence we have:

\begin{proposition}\label{Correspondence between monotone maps and modal intuitionistic frames}
    $\Box$-frames $\mathfrak{X}=(X,\leq,R)$, are in 1-1 correspondence with monotone maps  $R_{\mathfrak{X}}:X\to \mathsf{Up}(X)$ defined by
    \begin{equation*}
        x\mapsto R[x].
    \end{equation*}
\end{proposition}

This suggests that the category of coalgebras for the endofunctor $\mathsf{Up}(-)$ should be equivalent to the category of modal intuitionistic frames with modal p-morphisms. However, following the discussion in \cite{Litak2014}, note that we are working within the category $\mathbf{Pos}$, where maps are only required to be monotone, which means that there may be coalgebra morphisms which are not p-morphisms for the $\leq$-relation. If one tentatively restricts to $\mathbf{Pos}_{p}$, then the problem of coalgebra morphisms failing to be p-morphisms can be avoided. Indeed, $\mathsf{Up}(-)$ restricts to an endofunctor on this (non-full) subcategory. However, this means that the above assignment $R_{X}:X\to \mathsf{Up}(X)$ might no longer be a morphism in this category, since there is no guarantee that it will be a p-morphism. 

Hence, despite the intuitive  connection between the functor $\mathsf{Up}(-)$ and the semantics of intuitionistic modal logic, something seems to be missing for a coalgebraic representation. We will return to this in the next section. For now we will need to also discuss how this plays out in the context of general frames.

\subsection{Descriptive general frames for IML}

Just like in classical modal logic, in order to address the phenomenon of Kripke incompleteness, more general structures are required, in the form of \textit{intuitionistic general frames}. We recall here this semantics \cite{Wolter1999}.

\begin{definition}
Let $(X,\leq,R)$ be a triple where $(X,\leq)$ is an Esakia space and $R\subseteq X\times X$. We say that $(X,\leq,R)$ is a \textit{modal Esakia space} if:
\begin{enumerate}
    \item Whenever $U$ is a clopen upset, then $\Box_{R}U$ is a clopen upset, where $\Box_{R}U=\{x\in X : R[x]\subseteq U\}$.
    \item For each $x\in X$, $R[x]$ is a closed upset.
\end{enumerate}
We call the quadruple $(X,\leq,R,\mathcal{A})$ where $(X,\leq,R)$ is a modal Esakia space and $\mathcal{A}=\mathsf{ClopUp}(X)$ a \textit{descriptive intuitionistic modal frame} (or $\Box$-general frame, for short).

Given a map $f:X\to Y$ between $\Box$-general frames, we say that this is a modal p-morphism if it is a continuous p-morphism with respect to both relations. We denote by $\mathbf{DiG}$ the category of $\Box$-general frames with modal p-morphisms.
\end{definition}

We will need the following fact, which is derived essentially from the persistence condition on valuations:

\begin{lemma}\label{Mix law lemma}
    If $(X,\leq,R,\mathcal{A})$ is a $\Box$-general frame, then $R$ and $\leq$ satisfy the following mix law:
    \begin{equation*}
        {R}{=}{\leq} \circ {R} \circ {\leq}
    \end{equation*}
\end{lemma}

Consequently, we can think of a $\Box$-general frame as a modal intuitionistic frame equipped with a compatible topology, or as an Esakia space equipped with a compatible relation $R$. In the latter case we say that $(X,R,\leq,\mathcal{A})$ is a $\Box$-general frame over $X$.

By combining the classical J\'onnson-Tarski duality and Esakia duality, $\Box$-general frames have been shown to provide a general completeness result for $\mathsf{IPC}_{\Box}$ with respect to $\Box$-general frames (see e.g. \cite{Palmigiano2004,Wolter1999,Celani1997}). 

We consider the following variation of the previously outlined Vietoris endofunctor: given an Esakia space $X$, we write
\begin{equation*}
    \mathsf{V}^{\uparrow}(X)\coloneqq \{C\subseteq X : C \text{ is a closed upset}\}
\end{equation*}
with a topology given by a subbasis of sets of the form
\begin{equation*}
    [U],\langle X-V\rangle
\end{equation*}
where $U,V$ range over clopen upsets. Then we have the following fact, a proof of which can be found in the Appendix\footnote{We include this proof for completeness. It actually follows from the results mentioned, for e.g. in \cite{Bezhanishvili2023remarksonhyperspaces}, where $(V^{\uparrow}(X),\subseteq)$ is noted to be a Priestley space, by realising that the order-dual of a Priestley space is again a Priestley space.}.

\begin{proposition}
    Given a Priestley space $X$, the order-topological space $(V^{\uparrow}(X),\supseteq)$ is again a Priestley space.
\end{proposition}
\begin{proof}
    Note that this satisfies the Priestley separation axiom, since the original space does: if $C\npreceq D$, then $C\nsupseteq D$. Hence there is clopen upset such that $D\subseteq U$ and $C\nsubseteq U$; this means that $C\in [U]$ an $D\notin [U]$, which was to show.

    Now we show compactness. Assume that
    \begin{equation*}
        \mathsf{V}^{\uparrow}(X)=\bigcup_{i\in I}[U_{i}]\cup \bigcup_{j\in J}\langle X-V_{j}\rangle
    \end{equation*}
    where $U_{i},V_{j}$ are clopen upsets. Look at $C=X-\bigcup_{j\in J}X-V_{j}$. Now, if $C=\emptyset$, then $X$ is covered by $X-V_{j}$, so we can extract a finite subcover, say $X=X-V_{0}\cup...\cup X-V_{n}$. Now if $A\in \mathsf{V}^{\uparrow}(X)$, then $A$ must intersect one of the above subsets, since it is non-empty, and so
    \begin{equation*}
        \mathsf{V}^{\uparrow}(X)=\langle X-V_{0}\rangle \cup...\cup \langle X-V_{n}\rangle.
    \end{equation*}
    Otherwise $C\neq \emptyset$, so $C\in \mathsf{V}^{\uparrow}(X)$. By construction, $C\in [U_{i}]$ for some $i\in I$, so because it is closed, using compactness we have obtain a finite subcover of $X-U_{i}$,
    \begin{equation*}
        X-U_{i}\subseteq X-V_{0}\cup...\cup X-V_{n}.
    \end{equation*}
    And then we can show that
    \begin{equation*}
        \mathsf{V}^{\uparrow}(X)=[U_{i}]\cup \langle X-V_{0}\rangle\cup...\cup \langle X-V_{n}\rangle,
    \end{equation*}
    which again shows compactness.
\end{proof}

In analogy with what we noted for intuitionistic modal frames, we can at this point note the following (see e.g. \cite{deGroot2020}):

\begin{proposition}\label{DIM-frames to Priestley morphisms}
$\Box$-general frames $\mathfrak{X}=(X,\leq,R,\mathcal{A})$ are in 1-1 correspondence with Priestley morphisms $R_{\mathfrak{X}}:X\to \mathsf{V}^{\uparrow}(X)$ defined by
\begin{equation*}
    x\mapsto R[x].
\end{equation*}
\end{proposition}
\begin{proof}
    First note that if $x\leq y$ then from $yRz$, by Lemma \ref{Mix law lemma}, we have $xRz$. So we have that $R[x]\supseteq R[y]$, which means that the map is monotone. To see that it is continuous, note that if $U\subseteq X$ is a clopen upset, then
    \begin{equation*}
        R_{\mathfrak{X}}^{-1}[[U]]=\{x : R[x]\in [U]\}=\{x : R[x]\subseteq U\}=\Box_{R}U
    \end{equation*}
    which is a clopen upset by assumption.

    Conversely, assume that $f:X\to \mathsf{V}^{\uparrow}(X)$ is a Priestley morphism. Then define a relation $R$ as follows:
    \begin{equation*}
        xRy \iff y\in f(x).
    \end{equation*}
    Note that by construction  $R[x]=f(x)$ will be a closed upset. Moreover, if $U$ is a clopen upset, then we have
    \begin{equation*}
        f^{-1}[[U]]=\{x : f(x)\in [U]\}=\Box_{R}U
    \end{equation*}
    is clopen, since $f$ is continuous. It is clear that these two assignments are each other's inverses.
\end{proof}

Just like before, one can ask whether $\mathsf{CoAlg}(V^{\uparrow})$ is equivalent to the category of $\Box$-frames, and the exact same pattern repeats here: not all Priestley morphisms will be p-morpshims, so one may wish to restrict to $\mathbf{Esa}$; and whilst $\mathsf{V}^{\uparrow}(-)$ is an endofunctor on Esakia spaces, the assignment $R_{X}:X\to \mathsf{V}^{\uparrow}(X)$ may not be a p-morphism.

Having gotten to this point, it is natural to wonder if indeed one needs different kinds of frames, possibly induced by a different kind of functor. As we will show, however, $\Box$-frames are enough, if one considers the right endofunctor on the category of Esakia spaces. This will be the subject of the next section.

\section{Generalizing Ghilardi}\label{Generalizing Ghilardi}

In this section we recover the key technical tools which will be needed in the sequel. These results can be found in \cite{almeidanoteoncolimitsofHeytingalgebras} and generalize the ideas from Ghilardi's classical construction of the free Heyting algebra through a step-by-step method; see also \cite{Bezhanishvili2011heytingcoalgebra,Gool2011ConstructingTL} for more details on this perspective. Here, reversing the order of previous sections, we start by handling the topological case:

\subsection{The functor $V_{G}$}

\begin{definition}
Let $X,Y,Z$ be Priestley spaces, and $g:X\to Y$ and $f:Y\to Z$ be Priestley morphisms. We say that $f$ is \textit{open relative to $g$}\footnote{We recall that p-morphisms between posets $X$ and $Y$ correspond precisely to the open maps between the underlying topological spaces of $X,Y$ when given the Alexandroff topology. The terminology here, standard in the literature, derives from this analogy.} ($g$-open for short) if it satisfies the following:
\begin{equation}
\tag{*} \forall a\in X, \forall b\in Y, (f(a)\leq b \implies \exists a'\in X, (a\leq a' \ \& \ g(f(a'))=g(b)).
\end{equation}
Given $S\subseteq X$, we say that $S$ is \textit{rooted} if there is a point $x\in S$ such that for each $y\in S$, we have $x\leq y$. We say that $S\subseteq X$, a closed subset, is \textit{g-open} (understood as a poset with the restricted partial order relation) if the inclusion is itself $g$-open. Equivalently:
\begin{equation*}
\forall s\in S, \forall b\in X (s\leq b \implies \exists s'\in S (s\leq s' \ \& \ g(s')=g(b)).
\end{equation*}
\end{definition}

\begin{definition}
    Let $g:X\to Y$ be a map between Priestley spaces. Then consider
    \begin{equation*}
        V_{g}(X)\coloneqq \{C\subseteq X : C \text{ is closed, rooted and $g$-open }\},
    \end{equation*}
    with the topology given by a subbasis consisting of sets of the form
    \begin{equation*}
        [U],\langle V\rangle
    \end{equation*}
    where $U,V$ are clopen subsets of $X$.
\end{definition}

The following is proven in \cite[Lemmas 10 and 11]{almeidanoteoncolimitsofHeytingalgebras}:

\begin{proposition}
    Given $g:X\to Y$ a Priestley morphism, the order-topological space $(V_{g}(X),\supseteq)$ is a Priestley space, equipped with a Priestley surjection $r_{g}:V_{g}(X)\to X$ sending each rooted subset to its root.
\end{proposition}

We refer to $r_{g}$ as the \textit{root map} with respect to $V_{g}(X)$. The key property which this construction enjoys which we will need is the following:

\begin{lemma}\label{Duality Lemma for Key Property}
Given a Priestley morphism $g:X\to Y$, and given a Priestley space $Z$ with a $g$-open Priestley morphism $h:Z\to X$, there exists a unique $r_{g}$-open, continuous and monotone map $h'$ such that the triangle in Figure \ref{fig:commutingtriangleofpriestleyspaces} commutes.

\begin{figure}[h]
\centering
\begin{tikzcd}
Z \arrow[rd, "h"'] \arrow[rr, "h'"] &   & V_{g}(X) \arrow[ld, "r"] \\                    & X &                         
\end{tikzcd}
\caption{Commuting Triangle of Priestley spaces}
\label{fig:commutingtriangleofpriestleyspaces}
\end{figure}
\end{lemma}

\begin{definition}
    Let $g:X\to Y$ be a Priestley morphism. The \textit{$g$-Vietoris complex}  $(V_{\bullet}^{g}(X),\leq_{\bullet})$ over $X$, is a sequence
    \begin{equation*}
        (V_{0}(X),V_{1}(X),...,V_{n}(X),...)
    \end{equation*}
    connected by morphisms $r_{i}:V_{i+1}(X)\to V_{i}(X)$ such that:
    \begin{enumerate}
        \item $V_{0}(X)=Y$ and $V_{1}(X)=X$;
        \item $r_{0}=g$;
        \item For $i>1$, $V_{i+1}(X)\coloneqq V_{r_{i}}(V_{i}(X))$;
        \item For $i>0$ $r_{i+1}=r_{r_{i}}:V_{i+1}(X)\to V_{i}(X)$ is the root map.
    \end{enumerate}
    We denote the projective limit of this family (in the category $\mathbf{Pries}$) by $V_{G}^{g}(X)$ (by duality, and the fact that the category of distributive lattices is cocomplete, e.g. \cite{completenessofvarieties}). When $g$ is the terminal map to the one element poset, we often omit it.
\end{definition}

The proof of the following proposition, which will be instrumental in our work, can be found in the Appendix:

\begin{proposition}\label{Lifting of maps to Esakia morphisms}
    Let $X$ be an Esakia space, $Y$ a Priestley space, and assume that $f:X\to Y$ is a Priestley morphism. Then there is a unique Esakia morphism $\overline{f}:X\to V_{G}(Y)$, extending $f$. This is given as follows: the family $f_{n}:X\to V_{n}(Y)$, given by
    \begin{enumerate}
        \item $f_{0}=f$;
        \item $f_{n+1}(x)=f_{n}[{\uparrow}x]$;
    \end{enumerate}
    consisting of continuous functions, and $\overline{f}:X\to V_{G}(Y)$ is given by
    \begin{equation*}
        \overline{f}(x)=(f_{0}(x),f_{1}(x),...).
    \end{equation*}
\end{proposition}
\begin{proof}
Using Lemma \ref{Duality Lemma for Key Property} repeatedly, starting with the terminal map and proceeding along the roots, we get, a sequence defined by
\begin{enumerate}
    \item For every $x\in X$, $f_{0}(x)=f(x)$;
    \item For every $x\in X$, $f_{n+1}(x)=f_{n}[{\uparrow}x]$.
\end{enumerate}
By uniqueness of inverse limits, an extension $f^{\infty}:X\to V_{G}(X)$ given by
\begin{equation*}
    x \mapsto (f_{0}(x),f_{1}(x),f_{2}(x),...).
\end{equation*}
such that this map commutes with all the root maps. Such a map is certainly unique, so if we can show that it is an Esakia morphism, we are done. So assume that $f^{\infty}(x)\leq y$. Consider the following:
\begin{equation*}
    S={\uparrow}x\cup \{f_{n}^{-1}[y(n)] : n\in \omega\}.
\end{equation*}
This is a family of closed subsets in $X$, which is an Esakia space. Moreover, it has the finite intersection property: if we consider only finitely many elements, note that since $y(n)\subseteq f_{n}(x)$, then there is some $k$ such that $x\leq k$, and $f^{\infty}(k)$ agrees with $y$ up to the level $n$, i.e., the finite intersection ${\uparrow}x\cap f_{0}^{-1}[y(0)]\cap...\cap f_{n}^{-1}[y(n)]$ is non-empty. By compactness, there exists some $x'\in \bigcap S$, which means precisely that there is some $x\leq x'$ such that $f^{\infty}(x')=y$. This shows that $f^{\infty}$ is a p-morphism, as desired.
\end{proof}

Using this one obtains:

\begin{theorem}\label{Existence of right adjoint map}
    The assignment $V_{G}$ is an endofunctor on the category $\mathbf{Pries}$ of Priestley spaces and Priestley morphisms; indeed it is the right adjoint to the inclusion of the category $\mathbf{Esa}$ of Esakia spaces and Esakia morphisms into $\mathbf{Pries}$. 
\end{theorem}

\subsection{The functor $P_{G}$}\label{The functor PG}

Having the analysis of the previous section, we can ask whether something similar can be done for posets in general. Indeed this is the case, if one restricts to \textit{image-finite} posets; we will present the construction, and then comment a bit on this restriction.

\begin{definition}
    Let $g:X\to Y$ be a monotone map between posets. Then consider:
    \begin{equation*}
        P_{g}(X)=\{C\subseteq X : C \text{ is finite, rooted and $g$-open} \}.
    \end{equation*}
\end{definition}

\begin{definition}
    Let $g:X\to Y$ be a monotone map between image-finite posets. The $g$-discrete complex $(P_{\bullet}^{g}(X),\leq_{\bullet})$ over $X$, is a sequence
    \begin{equation*}
        (P_{0}(X),P_{1}(X),...,P_{n}(X),...)
    \end{equation*}
    connected by morphisms $r_{i}:P_{i+1}(X)\to P_{i}(X)$ such that
    \begin{enumerate}
        \item $P_{0}(X)=Y$ and $P_{1}(X)=X$;
        \item $r_{0}=g$;
        \item For $i>1$, $P_{i+1}(X)\coloneqq P_{r_{i}}(P_{i}(X))$;
        \item For $i>1$, $r_{i+1}\coloneqq r_{r_{i}}:P_{i+1}(X)\to P_{i}(X)$ is the root map.
    \end{enumerate}
    We denote the \textit{image-finite part} of the projective limit of this family (in $\mathbf{Pos}$) by $P^{g}_{G}(X)$. When $g$ is the terminal map to the one element poset we often omit it.
\end{definition}

The following propositions are the analogues of Propositions \ref{Lifting of maps to Esakia morphisms} and Theorem \ref{Existence of right adjoint map}.

\begin{proposition}\label{Lifting of maps to image-finite posets}
    Let $X$ be an image-finite poset, $Y$ a poset, and assume that $f:X\to Y$ is a monotone map. Then there is a unique p-morphism $\overline{f}:X\to P_{G}(Y)$, extending $f$. This is given as follows: the family $f_{n}:X\to P_{n}(Y)$, given by
    \begin{enumerate}
        \item $f_{0}=f$;
        \item $f_{n+1}(x)=f_{n}[{\uparrow}x]$;
    \end{enumerate}
    consists of monotone maps, and $\overline{f}:X\to V_{G}(Y)$ is given by
    \begin{equation*}
        \overline{f}(x)=(f_{0}(x),f_{1}(x),...)
    \end{equation*}
\end{proposition}

\begin{theorem}\label{Existence of right adjoint map for discrete}
    The assignment $P_{G}$ is an endofunctor on the category $\mathbf{Pos}$ of posets and monotone maps; indeed it is the right adjoint to the inclusion of the category $\mathbf{ImFinPos}_{p}$ of Image-finite posets spaces and p-morphisms into $\mathbf{Pos}$. 
\end{theorem}

\begin{remark}
    The reader may wonder why the restriction is done to the image-finite case. For the purpose of the results of this section, the issue lies in proving Proposition \ref{Lifting of maps to image-finite posets}: without the restriction to image-finiteness, it is not clear that the lifting $\overline{f}(x)$ will indeed be a p-morphism. But there are more general considerations which make this a natural restriction: as discussed in \cite[Section 7]{almeidanoteoncolimitsofHeytingalgebras}, there are several category-theoretic facts which make such a category much better behaved -- for instance, it is monadic over the category of posets, precisely through the above construction, and it corresponds to the profinite completion of the category of finite posets with p-morphisms. By contrast, the category $\mathbf{Pos}_{p}$ does not immediately seem to enjoy such properties.
\end{remark}

In light of this remark, the restriction to $\mathbf{ImFinPos}_{p}$ will be assumed throughout this paper, and the question of how to provide coalgebraic representations for arbitrary posets is left open.

\section{Main Results}\label{Main Results}

In this section we show that the category $\mathbf{CoAlg}(\mathsf{V}_{G}(\mathsf{V}^{\uparrow}(-)))$ is equivalent to the category $\mathbf{DiG}$ of descriptive $\Box$-frames with modal p-morphisms. In particular,  we will prove the following:

\begin{theorem}
Let $(X,\leq)$ be an Esakia space. Then the following are in one-to-one correspondence:
\begin{enumerate}
    \item $\Box$-frames over $X$;
    \item Priestley morphisms $f:X\to \mathsf{V}^{\uparrow}(X)$;
    \item Esakia morphisms $f':X\to \mathsf{V}_{G}(\mathsf{V}^{\uparrow}(X))$.
\end{enumerate}
\end{theorem}
\begin{proof}
    (i) is equivalent to (ii) by Lemma \ref{DIM-frames to Priestley morphisms}. For (iii) to (ii), if $f':X\to \mathsf{V}_{G}(\mathsf{V}^{\uparrow}(X))$ is an Esakia morphism, then certainly the map
    \begin{equation*}
        f\coloneqq f'\circ \pi_{0}
    \end{equation*}
    which projects everything to the first coordinate is a Priestley morphism. So we focus on the implication from (ii) to (iii). Given the map $f:X\to \mathsf{V}^{\uparrow}(X)$, by
    Proposition \ref{Lifting of maps to Esakia morphisms} there is a unique Esakia morphism $\overline{f}:X\to \mathsf{V}_{G}(\mathsf{V}^{\uparrow}(X))$ extending it. It is moreover clear, by definition, that these two assignments are each other's inverses.
\end{proof}

We moreover have the following:

\begin{theorem}
    The category $\mathbf{DiG}$ is equivalent to the category $\mathbf{CoAlg}(\mathsf{V}_{G}(\mathsf{V}^{\uparrow}(-)))$.
\end{theorem}
\begin{proof}
    Similarly to above, we show that if $f:X\to Y$ is a modal p-morphism, then we can lift this to a coalgebra morphism making the right diagram commute. Let $i_{X}:X\to \mathsf{V}_{G}(\mathsf{V}^{\uparrow}(X))$ and $i_{Y}:Y\to \mathsf{V}_{G}(\mathsf{V}^{\uparrow}(Y))$ be the coalgebra maps. It is clear, since $\mathsf{V}_{G}(\mathsf{V}^{\uparrow}(-))$ map is a functor, that this lifts to a map, which essentially depends on the map
    \begin{equation*}
        f^{*}:\mathsf{V}^{\uparrow}(X)\to \mathsf{V}^{\uparrow}(Y).
    \end{equation*}
    Indeed, this amounts, for $x\in X$, to have that $f^{*}(\pi_{0}[i_{X}(x)])=f[R[x]]$ equal to $R[f(x)]$, since the liftings will commute with these maps; but it is clear that such an equality means precisely that $f$ is a modal p-morphism.
\end{proof}

Exactly the same arguments, using Proposition \ref{Correspondence between monotone maps and modal intuitionistic frames}, and using the appropriate discrete versions of the results shown in Section \ref{Generalizing Ghilardi}, provide the following:

\begin{theorem}
    The category $\mathbf{ImFinK}$ of image-finite modal Kripke frames is equivalent to the category $\mathbf{CoAlg}(P_{G}(\mathsf{Up}(-)))$.
\end{theorem}

\begin{remark}
    The above result can now explain the phenomenon we alluded to in Section \ref{Preliminaries}. The noted correspondence between $\Box$-frames and coalgebras for $\mathsf{V}^{\uparrow}$, witnesses an equivalence between the category $\mathbf{DiG}^{m}$ of $\Box$-general frames with Priestley morphisms satisfying the p-morphism conditions for $R$, and the category of coalgebras for that functor. This provides a coalgebraic representation for positive modal logic over $\Box$. But this should not suffice to model the implication. One way to do so, is to introduce a construction which freely adds the implications. This is exactly the role played by our functor $\mathsf{V}_{G}(\mathsf{V}^{\uparrow}(-))$.
\end{remark}

\section{Applications}\label{Applications}

In this section we provide a few applications of our characterizations, exploiting and expanding the coalgebraic ideas exposed above.

\subsection{Bisimulations of $\Box$-frames}

As an illustration of the correctness of this coalgebraic representation, we show that the notion of bisimulation one obtains corresponds to the ones we would naturally want for $\Box$-frames and $\Box$-general frames. We show this for the functor $P_{G}(\mathsf{Up}(-))$, though similar results could easily be derived for $V_{G}$, since the former is slightly simpler to handle.

\begin{definition}
    Let $\mathfrak{X}=(X,\leq,R)$ and $\mathfrak{Y}=(Y,\leq',R')$ be two $\Box$-frames. We say that a relation $\sim\subseteq X\times Y$ is a $\Box$-bisimulation if:
    \begin{enumerate}
        \item Whenever $x,y\in X$ and $xSy$ and $x\sim x'$ where $x'\in Y$, then there is some $y'\in Y$ such that $ySy'$ and $y\sim y'$;
        \item Whenever $x',y'\in Y$ and $x\in X$ and $x'Sy'$ and $x\sim x'$ then there is some $y\in X$ such that $y\sim y'$.
    \end{enumerate}
    Where $S$ ranges over $\leq$ and $R$.
\end{definition}

The following theorem, proven in the Appendix, then follows by similar techniques to those used above:

\begin{theorem}
    Let $(X,\leq,R)$ and $(Y,\leq,R)$ be two image-finite modal intuitionistic frames. Then the following are in one-to-one correspondence:
    \begin{enumerate}
        \item $\Box$-bisimulations between $X$ and $Y$;
        \item Bisimulations for the endofunctor $\mathsf{P}_{G}(\mathsf{Up}(-))$. 
    \end{enumerate}
\end{theorem}
\begin{proof}
    Recall that a bisimulation in the category of image-finite posets with p-morphisms, for the functor $\mathsf{P}_{G}(\mathsf{Up}(-))$ amounts to a relation $B\subseteq X\times Y$ endowed with a coalgebra structure making the following diagram commute:
\begin{figure}[h]
    \centering
\begin{tikzcd}
X \arrow[d]    & B \arrow[d] \arrow[r, "\pi_{Y}"] \arrow[l, "\pi_{X}"'] & Y \arrow[d]    \\
\mathsf{P}_{G}(\mathsf{Up}(X)) & \mathsf{P}_{G}(\mathsf{Up}(B)) \arrow[r] \arrow[l]                     & \mathsf{P}_{G}(\mathsf{Up}(Y))
\end{tikzcd}    \caption{Bisimulations for the $P_{G}$}
    \label{fig:bisimulationupsetfunctor}
\end{figure}

Now given the bisimulation $B$, first define a structure on $X\times Y$ by sending $(x,y)$ to $R[x]\times R[y]$; then lift this to a map to $P_{G}(\mathsf{Up}(X\times Y))$ using Proposition \ref{Lifting of maps to image-finite posets}; the uniqueness ensures that the desired map commutes.

Conversely, assume that we have $B\subseteq X\times Y$ yielding a commuting diagram. Then we claim that $B$ is a $\Box$-bisimulation. Because the projetions are required to be $p$-morphisms, the bisimulation satisfies the clauses for the $\leq$-relation, and the fact that it does so for the $R$-relation amounts to projecting onto the first coordinate.   
\end{proof}

\subsection{Constructing the Free Intuitionistic Modal algebra}

The step-by-step construction of free algebras has long been employed in constructing modal algebras. This amounts to giving access to these algebras in a way that exploits the simplicity of working with simpler algebras, like Boolean algebras. In addition to that, these methods typically depend on the finitarity of some underlying algebra, and it is not always clear how to extend them to the infinite case (see \cite{Bezhanishvili2011heytingcoalgebra} for some discussion of this). Algebraically, the mechanism at play is the existence of a locally finite reduct over which one can ``layer" a non-locally finite operation. However, when in the face of two non-locally finite operations -- such as when having Heyting implications and modalities -- the situation can quickly become difficult. This is what we will be concerned with in this section: the construction of free modal Heyting algebras.

\begin{definition}
    Let $X$ be an Esakia space. Define the following sequence:
    \begin{equation*}
        (M_{0}(X),M_{1}(X),...,M_{n}(X),...)
    \end{equation*}
    and a sequence of morphisms $\pi_{k}:M_{k}(X)\to M_{k-1}(X)$, for $k>0$ and $\pi_{0}:M_{0}(X)\to M_{0}(X)$ as follows:
    \begin{enumerate}
        \item $M_{0}(X)=X$;
        \item $M_{n+1}(X)\coloneqq X\times \mathsf{V}_{G}(\mathsf{V}^{\uparrow}(M_{n}(X)))$;
        \item $\pi_{0}=id_{M_{0}}$;
        \item $\pi_{n+1}(x,C)=(x,\pi_{n}[C])$.
    \end{enumerate}
For each $k$, define the relation $R_{k}\subseteq (X\times \mathsf{V}_{G}(\mathsf{V}^{\uparrow}(M_{k}(X))\times M_{k}(X)$, given by
\begin{equation*}
    (x,C)R_{k}y \iff y\in \pi_{0}[C].
\end{equation*}
    Let $M_{\infty}(X)$ be the inverse limit (in the category of Priestley spaces) of these spaces. Define a relation $R_{\omega}$ as follows: if $x,y\in M_{\infty}(X)$
    \begin{equation*}
        xR_{\omega}y \iff \forall k\in \omega, x(k+1) R_{k} y(k).
    \end{equation*}
\end{definition}

First we note the following:

\begin{proposition}
    Given any Esakia space $X$, $M_{\infty}(X)$ is a modal Esakia space.
\end{proposition}
\begin{proof}
    The fact that the inverse limit is an Esakia space follows straightforwardly from duality. Moreover, given any point $x\in M_{\infty}(X)$ $R_{\omega}[x]=\bigcap\{y : x(k+1)R_{k}y(k)\}$, so to show this is point-closed it suffices to show that $R_{k}$ is point-closed. But this amounts to the fact that $x(k+1)=(z,C)$ where $C\in \mathsf{V}_{G}(\mathsf{V}^{\uparrow}(M_{k}(X)))$, and so $y\in R_{k}[x(k+1)]$ if and only if $y\in \pi_{0}[C]$, i.e., $R_{k}[x(k+1)]=\pi_{0}[C]$, which is closed by construction. Similarly, if $U\subseteq M_{\infty}(X)$ is a clopen downset, then note that by construction, for some $n$, $\pi_{n}[U]\subseteq M_{n}(X)$ is a clopen downset. Then consider    
    \begin{equation*}
        \Box_{R_{n+1}}[\pi_{n}[U]]=\{x\in M_{n+1}(X) : \forall y\in M_{n}( xR_{n+1}y \rightarrow y\in \pi_{n}[U])\}.
    \end{equation*}
    Since $xR_{n+1}y$ holds if $y\in \pi_{0}[C]$, this amounts, after unfolding the definitions, to
    \begin{equation*}
        [\pi_{n}[U]]=\{C\subseteq \mathsf{V}^{\uparrow}(M_{n}(X)) : C\subseteq \pi_{n}[U]\},
    \end{equation*}
    which is clopen upset by construction. By the topology of inverse limits, this implies that $\Box_{R_{\omega}}U$ is clopen upset as well.
\end{proof}

Now suppose that $Y$ is a modal Esakia space, and assume that $p:Y\to X$ is a monotone map. Then define the map $p_{1}:Y\to \mathsf{V}^{\uparrow}(X)$ as follows:
\begin{equation*}
    \overline{p}(y)=p[R[y]];
\end{equation*}
this is well defined, since if $z\in R[y]$, and $p(z)\leq w$, then by the p-morphism condition, there is some $z'$ such that $p(z')=w$, and $z\leq z'$; then $z'\in R[y]$, as well, so $p[R[y]]$ is an upwards closed set, and closed as well. To see that this is continuous, note that if $U$ is a clopen upset, then
\begin{equation*}
    \overline{p}^{-1}[[U]]=\{y\in Y : p[R[y]]\subseteq U\}=\{y\in Y : R[y]\subseteq p^{-1}[U]\}=\Box p^{-1}[U]
\end{equation*}
and since $p^{-1}[U]$ is a clopen upset, the latter is clopen because $Y$ is a modal Esakia space. Using Proposition \ref{Lifting of maps to Esakia morphisms}, we obtain a p-morphism $\overline{\overline{p}}:Y\to V_{G}(\mathsf{V}^{\uparrow}(X))$, and hence a map
\begin{align*}
    p_{1}:Y &\to X\times \mathsf{V}_{G}(\mathsf{V}^{\uparrow}(X))\\
    y &\mapsto (p(y),\overline{\overline{p}}(y)).
\end{align*}
We thus define a sequence of p-morphisms $p_{n}:Y\to M_{n}(X)$, which in turn induces a unique map $p_{\infty}:Y\to M_{\infty}(X)$, which is likewise a p-morphism. Then we can show the following:

\begin{proposition}
    Given $X$ an Esakia space and $Y$ a modal Esakia space, and $p:Y\to X$ a p-morphism, the unique lifting $p_{\infty}:Y\to M_{\infty}(X)$ is a modal Esakia morphism.
\end{proposition}
\begin{proof}
It suffices to show that if $xRy$, then $p_{k+1}(x)R_{k}p_{k}(y)$. To see this, in turn it suffices to show that $p_{k}(y)\in \pi_{0}[p_{k+1}(x)]$, i.e., $p_{k}(y)\in \overline{p_{k}}(x)$. But this is given by definition, since $\overline{p_{k}}(x)=p_{k}[R[x]]$. Hence we have the result.
\end{proof}

\begin{theorem}
    Given $X$ a finite set\footnote{The restriction to finiteness is insubstantial; if one wishes to consider arbitrary sets, replace $\mathcal{P}(X)$ by the Priestley space dual to the free distributive lattice generated by $X$ generators.}, consider $\mathcal{P}(X)$ as a poset with reverse inclusion. Then  $M_{\infty}(\mathcal{P}(X))$ is the dual to the free modal Heyting algebra on $X$ many generators.
\end{theorem}

\subsection{Coalgebraic Intuitionistic Logic}

The results of the previous sections point to a possible avenue for coalgebraic intuitionistic logic. To illustrate this we will provide here a coalgebraic semantics for intuitionistic neighbourhood frames. This is done in the case of image-finite posets, to maintain consistency with the available results in the literature, and for simplicity.

Recall that the logic $\mathsf{IPC}_{\Box,N}$ is defined over the same language as $\mathsf{IPC}$, but omitting the two normality axioms. An \textit{intuitionistic neighbourhood frame} is a triple $(X,\leq,N)$ of a poset together with a monotone map $N:X\to \mathcal{P}(\mathsf{Up}(X))$, where $\mathcal{P}(-)$ is ordered by inclusion (see e.g. \cite{Groot2020GoldblattThomasonTF} for a related approach, see also \cite{Dalmonte2020}). The morphisms between such frames are functions $f:X\to X'$ satisfying
\begin{equation*}
    a'\in N'(f(x)) \iff f^{-1}(a')\in N(x)
\end{equation*}
for all $x\in X$ and $a'\subseteq X'$. We denote by $\mathbf{ImFinN}$ the category of image-finite neighbourhood $\Box$-frames.

Consider the assignment:
\begin{equation*}
    \mathcal{P}(\mathsf{Up}(-)):\mathbf{Pos}\to \mathbf{Pos},
\end{equation*}
which sends a poset $P$ to the powerset (with inclusion) of its upset (with reverse inclusion), and sends monotone maps to their direct image. Then this is in an endofunctor on $\mathbf{Pos}$. Just like before, we can consider the composite functor $\mathsf{P}_{G}(\mathcal{P}(\mathsf{Up}(-)))$. Then we have:

\begin{theorem}
    There is an equivalence between $\mathbf{CoAlg}(\mathsf{P}_{G}(\mathcal{P}(-)))$ and the category $\mathbf{ImFinN}$.
\end{theorem}

The proof of this theorem follows by using exactly the same tools as before, and would go through just as well if we replaced $\mathbf{ImFinN}$ by the corresponding category of neighbourhood $\Box$-general frames, and the functors by their appropriate Vietoris-style variations. In fact, it suggests a more general phenomenon which we now point out:

\begin{definition}
Let $F:\mathbf{Pries}\to \mathbf{Pries}$ be an endofunctor on the category of Priestley spaces; we define $F^{*}:\mathbf{Esa}\to \mathbf{Esa}$: the \textit{intuitionistic lifting} of $F$ to be the functor obtained by composition in the following diagram:
\begin{figure}[h]
    \centering
\begin{tikzcd}
\mathbf{Pries} \arrow[r, "F"]                          & \mathbf{Pries} \arrow[d, "P_{G}"] \\
\mathbf{Esa} \arrow[r, "F^{*}"'] \arrow[u, "I"] & \mathbf{Esa}              
\end{tikzcd}    \caption{Intuitionistic Lifting of functor F}
    \label{fig:intuitionisticreflection}
\end{figure}

Similarly, given an endofunctor $F:\mathbf{Pos}\to \mathbf{Pos}$ on the category of posets with monotone maps, we define $F^{*}:\mathbf{ImFinPos_{p}}\to \mathbf{ImFinPos_{p}}$ the \textit{intuitionistic reflection} of $F$ to be the functor obtained by composing the obvious diagram.
\end{definition}

The results presented so far indicate a way to move from \textit{positive distributive logics} to \textit{intuitionistic logics}, using the mechanism of intuitionistic lifting. Hence it opens an avenue into investigations of intuitionistic coalgebraic logic. We leave a systematic study of the properties of intuitionistic lifting of functors, in a coalgebraic setting, for further work.

\section{Conclusions and Future Work}\label{Conclusions}

In this paper we established a coalgebraic representation for descriptive intuitionistic modal frames, and for image-finite modal Kripke frames. Using this we provided a construction of the free modal Heyting algebra generated by an Esakia space, which has as a special case the free modal Heyting algebra generated by a set of generators. The work presented leaves several questions open, of a technical and conceptual nature.

Of a more technical nature, there are several facts to be clarified about $P_{g}$ and $V_{g}$ which would be of interest; for instance, it is not clear, assuming that $X$ is an Esakia space, that $V_{g}(X)$ remains an Esakia space. Such technical facts would be relevant for the study of normal forms in intuitionistic modal logic.

In terms of the scope of the approach, there are naturally several possible lines of development. Like intuitionistic logic, $\mathsf{S4}$ modal logic is axiomatised by axioms of rank greater than $1$, which has so far kept it from being represented coalgebraically; we expect that similar techniques to the ones exposed here should account for these and similar cases. 
%Extensions to other settings such as bi-intuitionistic logic would likewise be very interesting to explore, and should still be possible to explore within this framework, though they may present additional problems and challenges.

Finally, related to our final remarks, the question of how to develop coalgebraic intuitionistic logic remains open, though the approach exposed here shows one recipe: first develop positive distributive logic, and then lift it to an intuitionistic setting. 
%This naturally imposes some limitations on the interactions of implications and modalities which leave out several interesting calculi (for example, the calculus of Kuznetsov-Muravitsky, see e.g. \cite{Esakia2006} for an axiomatisation). 
We leave a full study of this situation for future research.

\section{Acknowledgements}

The authors would like to thank the two anonymous referees for suggestions that improved the presentation of this work.

%% Bibliography
%% Make sure to use the bibliographystyle aiml22.
\bibliographystyle{aiml22}
\bibliography{aiml22}

\begin{thebibliography}{10}
\expandafter\ifx\csname url\endcsname\relax
  \def\url#1{\texttt{#1}}\fi
\expandafter\ifx\csname urlprefix\endcsname\relax\def\urlprefix{URL }\fi
\newcommand{\enquote}[1]{``#1''}

\bibitem{completenessofvarieties}
Adámek, J. and J.~Rosický, \emph{Which categories are varieties? ((co)algebraic pearls)} (2021).
\newline\urlprefix\url{https://drops.dagstuhl.de/entities/document/10.4230/LIPIcs.CALCO.2021.6}

\bibitem{almeidanoteoncolimitsofHeytingalgebras}
Almeida, R.~N., \emph{Colimits of heyting algebras through esakia duality} (2024).
\newline\urlprefix\url{https://arxiv.org/abs/2402.08058}

\bibitem{Balan2015}
Balan, A., A.~Kurz and J.~Velebil, \emph{Positive fragments of coalgebraic logics}, Logical Methods in Computer Science \textbf{Volume 11, Issue 3} (2015).
\newline\urlprefix\url{http://dx.doi.org/10.2168/LMCS-11(3:18)2015}

\bibitem{Bezhanishvili2023remarksonhyperspaces}
Bezhanishvili, G., J.~Harding and P.~Morandi, \emph{Remarks on hyperspaces for priestley spaces}, Theoretical Computer Science \textbf{943} (2023), p.~187–202.
\newline\urlprefix\url{http://dx.doi.org/10.1016/j.tcs.2022.12.001}

\bibitem{Bezhanishvili2011heytingcoalgebra}
Bezhanishvili, N. and M.~Gehrke, \emph{Finitely generated free heyting algebras via birkhoff duality and coalgebra}, Logical Methods in Computer Science \textbf{Volume 7, Issue 2} (2011).
\newline\urlprefix\url{https://doi.org/10.2168/lmcs-7(2:9)2011}

\bibitem{BezhanishviliK07}
Bezhanishvili, N. and A.~Kurz, \emph{Free modal algebras: {A} coalgebraic perspective}, in: T.~Mossakowski, U.~Montanari and M.~Haveraaen, editors, \emph{Algebra and Coalgebra in Computer Science, Second International Conference, {CALCO} 2007, Bergen, Norway, August 20-24, 2007, Proceedings},  Lecture Notes in Computer Science  \textbf{4624} (2007), pp. 143--157.
\newline\urlprefix\url{https://doi.org/10.1007/978-3-540-73859-6\_10}

\bibitem{Blackburn2002-fd}
Blackburn, P., M.~de~Rijke and Y.~Venema, \enquote{Modal logic,} Number~53 in Cambridge tracts in theoretical computer science, Cambridge University Press, Cambridge, England, 2002.

\bibitem{Celani1997}
Celani, S. and R.~Jansana, \emph{A new semantics for positive modal logic}, Notre Dame Journal of Formal Logic \textbf{38} (1997).
\newline\urlprefix\url{http://dx.doi.org/10.1305/ndjfl/1039700693}

\bibitem{Chagrov1997-cr}
Chagrov, A. and M.~Zakharyaschev, \enquote{Modal Logic,} Oxford Logic Guides, Clarendon Press, Oxford, England, 1997.

\bibitem{Dalmonte2020}
Dalmonte, T., C.~Grellois and N.~Olivetti, \emph{Intuitionistic non-normal modal logics: A general framework}, Journal of Philosophical Logic \textbf{49} (2020), p.~833–882.
\newline\urlprefix\url{http://dx.doi.org/10.1007/s10992-019-09539-3}

\bibitem{Davey2002-lr}
Davey, B.~A. and H.~A. Priestley, \enquote{Introduction to Lattices and Order,} Cambridge University Press, Cambridge, 2002, 175--200 pp.

\bibitem{Groot2020GoldblattThomasonTF}
de~Groot, J., \emph{Goldblatt-thomason theorems for modal intuitionistic logics}, ArXiv \textbf{abs/2011.10221} (2020).
\newline\urlprefix\url{https://api.semanticscholar.org/CorpusID:227118992}

\bibitem{deGroot2020}
de~Groot, J. and D.~Pattinson, \emph{Modal intuitionistic logics as dialgebraic logics}, in: \emph{Proceedings of the 35th Annual ACM/IEEE Symposium on Logic in Computer Science}, LICS ’20 (2020), pp. 355--369.
\newline\urlprefix\url{http://dx.doi.org/10.1145/3373718.3394807}

\bibitem{esakiatopologicalkripkemodels}
Esakia, L., \emph{Topological kripke models}, Dokl. Akad. Nauk SSSR \textbf{214} (1974), pp.~298--301.

\bibitem{Esakia2006}
Esakia, L., \emph{The modalized heyting calculus: a conservative modal extension of the intuitionistic logic}, Journal of Applied Non-Classical Logics \textbf{16} (2006), p.~349–366.
\newline\urlprefix\url{http://dx.doi.org/10.3166/jancl.16.349-366}

\bibitem{Esakiach2019HeyAlg}
Esakia, L., \enquote{Heyting Algebras: Duality Theory,} Springer, 2019, english translation of the original 1985 book.

\bibitem{FischerServi1977}
Fischer~Servi, G., \emph{On modal logic with an intuitionistic base}, Studia Logica \textbf{36} (1977), p.~141–149.
\newline\urlprefix\url{http://dx.doi.org/10.1007/BF02121259}

\bibitem{ghilardifreeheyting}
Ghilardi, S., \emph{Free heyting algebras as bi-heyting algebras}, C. R. Math. Rep. Acad. Sci. Canada  (1992), pp.~240--244.

\bibitem{Kupke2003}
Kupke, C., A.~Kurz and Y.~Venema, \emph{Stone coalgebras}, Electronic Notes in Theoretical Computer Science \textbf{82} (2003), p.~170–190.
\newline\urlprefix\url{http://dx.doi.org/10.1016/S1571-0661(04)80638-8}

\bibitem{Litak2014}
Litak, T., \enquote{Constructive Modalities with Provability Smack,} Springer Netherlands, 2014 p. 187–216.
\newline\urlprefix\url{http://dx.doi.org/10.1007/978-94-017-8860-1_8}

\bibitem{Palmigiano2004}
Palmigiano, A., \emph{A coalgebraic view on positive modal logic}, Theoretical Computer Science \textbf{327} (2004), p.~175–195.
\newline\urlprefix\url{http://dx.doi.org/10.1016/j.tcs.2004.07.026}

\bibitem{Servi1980}
Servi, G.~F., \enquote{Semantics for a Class of Intuitionistic Modal Calculi,} Springer Netherlands, 1980 p. 59–72.
\newline\urlprefix\url{http://dx.doi.org/10.1007/978-94-009-8937-5_5}

\bibitem{simpson1994proof}
Simpson, A.~K., \enquote{The proof theory and semantics of intuitionistic modal logic,} Ph.D. thesis, University of Edinburgh (1994).

\bibitem{Gool2011ConstructingTL}
van Gool, S., \emph{Constructing the lindenbaum algebra for a logic step-by-step using duality (extended version)}, in: \emph{Future Directions of Logic, Proceedings of PhDs in Logic III}, 2011, pp. 1--15.
\newline\urlprefix\url{https://api.semanticscholar.org/CorpusID:14809947}

\bibitem{Venema2007}
Venema, Y., \enquote{6 Algebras and coalgebras,} Elsevier, 2007 p. 331–426.
\newline\urlprefix\url{http://dx.doi.org/10.1016/S1570-2464(07)80009-7}

\bibitem{Venema2014modallogicandvietoris}
Venema, Y. and J.~Vosmaer, \enquote{Modal Logic and the Vietoris Functor,} Springer Netherlands, 2014 p. 119–153.
\newline\urlprefix\url{http://dx.doi.org/10.1007/978-94-017-8860-1_6}

\bibitem{Wolter1999}
Wolter, F. and M.~Zakharyaschev, \enquote{Intuitionistic Modal Logic,} Springer Netherlands, 1999 p. 227–238.
\newline\urlprefix\url{http://dx.doi.org/10.1007/978-94-017-2109-7_17}

\end{thebibliography}

\end{document}